\def\N{{\bf N}}

\def\CC{{\rm\kern.24em
   \vrule width.02em
       height1.4ex depth-.05ex
   \kern-.26em C}}
\def\QQ{{\rm\kern.24em
   \vrule width.02em
       height1.4ex depth-.05ex
   \kern-.26em Q}}
\def\PP{{\rm I\kern-.25em P}}         \def\RR{{\rm I\kern-.25em R}}
\def\DD{{\rm I\kern-.25em D}}         \def\EE{{\rm I\kern-.25em E}}
\def\FF{{\rm I\kern-.25em F}}         \def\NN{{\rm I\kern-.23em N}}
\def\RRp{{\rm I\kern-.25em R}_{+}}
\def\IND{{\rm 1\kern-.25em I}}

\documentclass[11pt,a4paper]{article}
\usepackage{graphicx}
\begin{document}
\font\gothic = eufm10 \font\Bbb = msbm10

\begin{center}
\Large{\bf A CONTINUOUS ANALOGUE OF THE INVARIANCE PRINCIPLE AND ITS ALMOST  SURE VERSION}
\end{center}

\bigskip

\centerline{By ELENA PERMIAKOVA (Kazan)}

Chebotarev inst. of Mathematics and Mechanics, Kazan State
University

Universitetskaya 17, 420008 Kazan

e-mail epermiakova@mail.ru

\medskip

{\small {\bf Abstract.} We deal with random processes
obtained from a homogeneous random process with independent
increments by replacement of the time scale and by multiplication
by a norming constant.
We prove the convergence in distribution of
these processes to Wiener process in the Skorohod space endowed by the topology of
uniform convergence.
An integral type almost sure version of this limit theorem is obtained.}

{\bf 2000 AMS Mathematics Subject Classification.}
60F05 Central limit and other weak theorems, 60F15 Strong theorems.

\medskip
{\bf Key words and phrases:} functional limit theorem, almost sure limit theorem, process with independent stationary increments.

\medskip

 \centerline{\bf 1. Introduction}

 \medskip

The usual invariance principle asserts the convergence of the sequence of the random processes
$X_n(x)=\frac{1}{\sqrt{n}}\sum_{i=1}^{[nx]}\xi_i$, $x\in{\bf R^+}$,
as $n\to\infty$, to the Wiener process $W$, where $\xi_i$ are
independent identically distributed centered random variables with variance $1$.
In this paper we study approximations of the Wiener
process $W$ by the random processes
$$
X_t(x)=\frac{1}{\sqrt{t}}V(tx), \,\,\,\,x\in{\bf R^+}, \eqno(1)
$$
where $t>0$ is a parameter and $V$ is  a centered  homogeneous
random process with independent increments such that $V(0)=0$ and
${\bf E}(V(1))^2=\sigma^2$.
Then almost all sample paths of $X_t$
belong to the Skorohod space $D[0,1]$.
In Section 2 we prove that   $X_t$
converges to $\sigma W$, as $t\to\infty$, in distribution in $D[0,1]$.

Almost sure versions of functional limit theorems were studied in several papers.
Here we mention only Lacey and Philipp [1], Chuprunov and Fazekas
[2], Chuprunov and Fazekas [3].

In Section 3 we prove almost  sure versions of our limit theorem.
Let $(\Omega, \mbox{\gothic A}, {\bf P})$ be the
probability space on which the random process $V$ is defined.
We show that for certain sequence $(s_n)$, the sequence of measures
$$
Q_{n}(\omega)=Q_{n}[X_{s_n}(t)](\omega)=
\frac{1}{\ln n}\sum_{k=1}^{n}\frac{1}{k}\delta _{X_{s_k}(t)(\omega)}  \eqno(1+)
$$
converges weakly the to the distribution of $\sigma W$ in
$D[0,1]$ for almost all $\omega\in \Omega$.
Here and in the following $\delta_{x} $ denotes the measure of unit mass,
concentrated in the point   $x$.

Also we prove  integral type almost sure versions of our limit
theorem. In Chuprunov and Fazekas [4] a general integral type
almost sure limit theorem is presented. Then the general theorem
is applied to obtain almost sure versions of limit theorems for
semistable and max-semistable processes, moreover for processes
being in the domain of attraction of a stable law or being in the
domain of geometric partial attraction of a semistable or a
max-semistable law. We mention a simple consequence of the general
result of Chuprunov and Fazekas [4]. Let $V(t)$ be a process with
characteristic function (3). Let $f$ be a function with property
(B). Let $X_t$ be defined by (1). Then
$$
\frac{1}{\log(T)}\int_1^T\delta_{X_{f(t)}(1)(\omega)}\frac{1}{t}dt
$$
converges weakly,  as $T\to\infty$, to centered  gaussian
distribution with the variance $K(\infty)$.

In this paper we prove a functional version of this proposition.
We show that the measures
$$
Q_T^I=\frac{1}{D(T)}\int_1^T\delta_{X_{f(t)}(\omega)}d(t)dt
$$
converge weakly, as $T\to\infty$, to the distribution of $\sigma W$
 in $D[0,1]$ for almost all $\omega\in \Omega$.
The proof of this result is based on the criterion for integral type almost sure
version of a limit theorem which was obtained in Chuprunov and Fazekas [4].

\bigskip

\centerline{\bf 2. Functional limit theorems}

\medskip

We will denote by $\stackrel{d}{\to}$ the convergence in distribution,
by $\stackrel{w}{\to}$ the weak convergence of measures,
by $\mu_{\zeta}$ the distribution of the random element $\zeta$ and
by $\mbox{\gothic B}({\bf B})$ the $\sigma$-algebra of
the Borel subsets of the metric space ${\bf B}$.
In the paper we will denote by the same symbols the random process and the random element
corresponding to this random process.

Using the Kolmogorov representation (see [5], sect. 18) we can
assume that the characteristic function of the centered  homogeneous
random process $V(t)$ with independent increments is
$$
\phi_{V(t)}(x)={\bf
E}\left(e^{ixV(t)}\right)=e^{\omega(t,x,K(x))}= \eqno(3)
$$
$$
=\exp\left(t\left\{\int_{-\infty}^{+\infty}(e^{ixy}-1-ixy)\frac{1}{y^2}dK(y)\right\}\right)
,\,\,\,\, x\in{\bf R}.
$$
Here $K(y)$ is an bounded increasing function such that $K(-\infty)=0$.

We will consider the sequence of the random processes
$$
Y_n= X_{s_n}, \eqno(4)
$$
where
$$
s_n\to\infty,\,\,\,\,{\rm as}\,\,\,\, n\to\infty. \eqno(5)
$$

We will use the following preliminary result.

\medskip
{\bf Theorem 1.} {\it Let $Y_n$ be defined by (1) and (4) and assume that
$(s_n)$ satisfies (5).
Let $\sigma^2=K(\infty)$.
Then we have
$$
Y_n\stackrel{d}{\to} \sigma W, \,\,\,{\it as}\,\,\,n\to\infty,
$$
in $D[0,1]$ endowed by the topology of uniform convergence.}

{\sc Proof.} Let $0\le t_1<t_2<\infty$.
By the convergence criterionn in [5], sect. 19,
$$
Y_n(t_2)-Y_n(t_1)\stackrel{d}{\to}
\sigma(W(t_2)-W(t_1)),\,\,\,\,{\rm as}\,\,\, \ \ n\to\infty. \eqno(6)
$$
 Let $0\le t_0<t_1<\dots <t_k<\infty$.
Introduce the notation $\Delta Y_{ni}=Y_n(t_i)-Y_n(t_{i-1})$
and $\Delta W_{i}=W(t_i)-W(t_{i-1})$.
Since  $\Delta Y_{ni}$, $1\le i\le k$,
are independent random variables, from (6) we obtain
$$
(\Delta Y_{n1},\dots, \Delta Y_{nk})\stackrel{d}{\to} (\sigma
\Delta W_{1},\dots, \sigma \Delta W_{k}).
$$
Consequently, the finite dimensional distributions of $Y_n$ converge
to the finite dimensional distributions of $\sigma W$.
Also we have
$$
\limsup_{n\to\infty}{\bf
E}|Y_n(t_2)-Y_n(t_1)|^2=\sigma^2|t_2-t_1|. \eqno(7)
$$

But (7) together with the convergence of the finite dimensional
distributions gives the weak convergence of $Y_n$ to $\sigma W$
(see Billingsley [6], Theorem 15.6) in $D[0,1]$ with Skorohod's $J_1$-topology.
However, in our case the limit process is a continuous one.
Hence, (see Pollard [7], p. 137, and the discussion in
Billingsley [6], Sect. 18), the weak convergence in Skorohod's $J_1$-topology actually implies the weak convergence
in the uniform topology of $D[0,1]$.
The proof is complete.

\medskip
Using Theorem 1, we can prove the following theorem.

\smallskip
{\bf Theorem 2.} {\it Let $X_t$ be defined by (1). Then it holds that
$$
X_t\stackrel{d}{\to} \sigma W, \,\,\,{\it as}\,\,\,t\to\infty,
$$
in $D[0,1]$ endowed by the topology of uniform convergence.}

{\sc Proof.} Consider $D[0,1]$ in the topology of uniform convergence
and the space $M$ of distributions on $D[0,1]$ with the topology of convergence in distribution.
Then $M$ is a metric space and denote by $\rho_M$ a metric which defines this topology on $M$.
Then, by Theorem 1, $\rho_M(\mu_{X_{s_n}}, \mu_{\sigma W})\to 0$, as $s_n\to\infty$.
Therefore $\rho_M(\mu_{X_t}, \mu_{\sigma W})\to 0$, as $t\to\infty$.
The proof is complete.

\medskip
Let $\xi_i$, $i\in {\bf N}$, be independent identically
distributed random variables with the expectation $a$ and the
variance $\sigma^2$ and let  $\pi(t)$, $t\in {\bf R}^+$, be a
Poissonian process with the intensity 1, and let the family
$\xi_i$, $i\in {\bf N}$,be independent of $\pi(t)$, $t\in {\bf
R}^+$.  Then
$$
V(x)=\sum\limits_{i=1}^{\pi(x)}\xi_i - ax, \,\,\,\,x\in{\bf R}^+,
$$
is  a centered  homogeneous random process with independent
increments such that $V(0)=0$ and ${\bf E}(V(1))^2={\sigma^2 +
a^2}$

So from Theorem 2 we obtain the following corollary.

\smallskip
{\bf Corollary 1.} {\it Let $\xi_i$, $i\in {\bf N}$, be
independent identically distributed random variables with the
expectation $a$ and the variance $\sigma^2$.
Let
$$
X'_t(x)=\frac{\sum\limits_{i=1}^{\pi(tx)}\xi_i - atx}{\sqrt{t}},
\,\,\,\,x\in{\bf R}^+. \eqno(8)
$$
Then one has
$$
X'_t\stackrel{d}{\to} \sqrt{\sigma^2 + a^2} W, \,\,\,{\it
as}\,\,\,t\to\infty,
$$
in $D[0,1]$ with the topology of uniform convergence.}

\medskip
For $\xi_i=1$ from Corollary 1 we obtain the following.

\smallskip
{\bf Corollary 2.} {\it Let
$$
X^*_t(x)=\frac{\pi(tx) - tx}{\sqrt{t}}, \,\,\,\,x\in{\bf R}^+. \eqno(9)
$$
Then one has
$$
X^*_t\stackrel{d}{\to}  W, \,\,\,{\rm as}\,\,\,t\to\infty,
$$
in $D[0,1]$ with the topology of uniform convergence.}

\bigskip

\centerline{\bf 3. Almost sure versions of functional limit theorems}

\medskip

We will consider the sequence of measures defined by (1+) and connected with the
random processes $X_{n}(t)$.

For the sequence $s_n$ we will assume  the following property
$$
{\rm for}\,\,\,\,{\rm some}\,\,\,\, \beta>0,\,\,\,\
\frac{s_n}{n^{\beta}}\,\,\,\,\, {\rm is}\,\,\, {\rm an}\,\,\, {\rm
increasing}\,\,\,{\rm sequence.} \leqno(A)
$$

{\bf Theorem 3.} {\it Let (A) be valid. Then it holds that
$$
Q_{n}(\omega)\stackrel{w}{\to}\mu_{\sigma W}\,\,\, \  \mbox{\rm if }\,\,\, \ n\to \infty \ \ \ \ \ \ \eqno(10)
$$
for almost all $\omega \in \Omega.$ }

{\sc Proof. } Let $l<k$. Let
$$
 Y_{kl}(x)=\left\{ \begin{array}{lr}
0,                                         & 0\leq x < \frac{s_l}{s_k},\\
Y_{k}(x) - \frac{V(s_l)}{\sqrt{s_k}}, \ \  & \frac{s_l}{s_k}\leq x \leq  1. \\
\end{array}
\right.
$$

Then
 $Y_{kl}(x)$, $0\le x\le 1$, and  $Y_{l}(x)$, $0\le x\le 1$, are independent random processes.
 Let $\rho$ be the metric of $D[0,1]$.
Using the moment inequality from [8], sect. 5, we obtain
$$
{\bf E} \rho (Y_{k},Y_{kl})\le
{\bf E} \sup_{0\leq x \leq 1}|Y_{k} (x)-Y_{kl}(x)|\le
{\bf E} \sup_{0\le x  \le \frac{s_l}{s_k}}|Y_{k}(x)|+
{\bf E} \left|\frac{V(s_l)}{\sqrt{s_k}}\right|\le
$$
$$
\le 4{\bf E}\left|Y_{k}\left(\frac{s_l}{s_k}\right)\right|+ {\bf
E} \left|\frac{V(s_l)}{\sqrt{s_k}} \right| \le 5\frac{\sqrt{{\bf
E}(V(s_l))^2}}{\sqrt{s_k}}
 \le 5\sigma \sqrt{\frac{s_l}{s_k}}
 \le 5\sigma
\left(\frac{l}{k}\right)^{\beta/2}. \eqno(11)
$$

By Lemma 1 from [3], this implies (10).
The proof is complete.

\bigskip
For the function $f$ we will consider following the property
$$
{\rm for}\,\,\,\,{\rm some}\,\,\, \beta>0, \,\,\,\, \
\frac{f(x)}{x^{\beta}}\,\,\,\, {\rm is}\,\,\, {\rm an}\,\,\, {\rm
increasing}\,\,\,{\rm function.} \leqno(B)
$$
\smallskip

Now we will prove the integral type almost sure version of Theorem 2.
We will consider the random processes
$$
Y_{t}(x)=\frac{V(f(t)x)}{\sqrt{f(t)}},\,\,\,\, \ 0\le x\le 1.
$$

We will assume that

\medskip\noindent
{\it (C)} \ \ \  the function $d(s)$ is a decreasing such that
$\int_{k}^{k+1} d(s)ds \leq \log \sqrt\frac{k+1}{k}$ for all $k\in{\bf N}$ and
$\int_{1}^{\infty}d(s)ds =+\infty. $
Let  $D(S)= \int_{1}^{S}d(s)ds.$

{\bf Theorem 4.} {\it Let (B) and (C) be valid. Then  we have
$$
Q_{S}^{I}(\omega)=\frac{1}{D(S)}\int_{1}^{S}\delta_{Y_{s}(\omega)}d(s)ds
\stackrel{w}{\to} \mu_{\sigma W}, \  \,\,\, {\it as}\,\,\, \ S\to \infty,
\eqno(12)
$$
for almost all $\omega \in \Omega$.}

{\sc Proof.} Let $0<l<k$, $l,k \in \N$,  $k \leq t \le k+1$.
Introduce the notation

$$
 Y_{lkt}(s)=\left\{ \begin{array}{lr}
0,                                         &  0\leq s \leq \frac{f(l+1)}{f(t)}, \\
Y_{t}(s)-\frac{V(f(l))}{\sqrt{f(t)},} \ \  &  \frac{f(l+1)}{f(t)}\leq s \leq  1. \\
\end{array}\right.
$$
Then $\{Y_{lkt}(s): k\leq t \leq k+1\}$ and $\{Y_{t}(s): l \leq t \leq l+1 \}$ are independent families.
 Repeating the proof of (11), we obtain
$$
{\bf E}\rho (Y_{t}, Y_{lkt})\le {\bf E}\sup_{0\leq s\leq
1}|X_{t}(s)-X_{lkt}(s)|\leq 5\cdot
2^{\beta/2}\sigma\left(\frac{l}{k}\right)^{\beta/2}.
$$
By Corollary 2.1, from Chuprunov and Fazekas [4] this and Theorem 2 implies (12).
The proof is complete.

\bigskip
Theorem 4 and Corollary 1 (resp. Corollary 2) of
Theorem 2 imply the corollaries.

{\bf Corollary 3.} {\it Let the $X'$ be defined by (8) and let $f$ be
a function with the property (B).
 Then
$$
\frac{1}{\ln(S)}\int_{1}^{S}\delta_{X'_{f(s)}(\omega)}
\frac{1}{s}ds \stackrel{w}{\to} \mu_{\sqrt{\sigma^2 + a^2}
W},\,\,\, \ {\it as}\,\,\,\  S\to \infty,
$$
for almost all $\omega \in \Omega.$}

{\bf Corollary 5.} {\it Let the $X^*$ be defined by (9) and let $f$ be
a function with the property (B.)
 Then
$$
\frac{1}{\ln(S)}\int_{1}^{S}\delta_{X^*_{f(s)}(\omega)}d(s)ds
\stackrel{w}{\to}
\mu_{W},\,\,\,\ {\it  as}\,\,\, \  S\to \infty,
$$
for almost all $\omega \in \Omega.$}

{\bf Remark 1.} Corollary 1 of Theorem 2 is a functional limit
theorem for random sums.
So  Corollary 3 of Theorem 4 is a
integral type almost sure version of a functional limit theorem for
random sums.
(For limit theorems for random sums see  Korolev and Kruglov [9].)

\bigskip

 \centerline{\bf  References}

 \medskip

 [1] {\sc M. T. Lacey} and {\sc W. Philipp,} A note on almost sure central limit theorem,
{\it Statistics and Probability Letters}
 {\bf9}(2) (1990), 201--205.

 [2] {\sc A. Chuprunov} and {\sc I. Fazekas},
Almost sure versions of some analogues of the invariance principle,
{\it Publicationes Mathematicae, Debrecen}
{\bf 54}(3-4) (1999), 457--471.

[3] {\sc A. Chuprunov} and {\sc I. Fazekas}, Almost sure limit theorems for the Pearson statistic,
{\it Teor. Veroyatnost. i Primenen.} {\bf 48}(1), 162--169.

[4] {\sc A. Chuprunov} and {\sc I. Fazekas},
Integral analogues of almost sure limit theorems,
{\it Periodica Mathematica Hungarica}
(to appear).

[5] {\sc B. V. Gnedenko} and {\sc A. N. Kolmogorov},
Limit Distributions for Sums of Independent Random Variables,
{\it  Addison-Wesley, Reading, Massachusetts},
1954.

[6] {\sc P. Billingsley}, Convergence of probability measures,
{\it John Wiley and Sons, New York}, 1968.

[7] {\sc D. Pollard}, Convergence of stochastic processes,
{\it Springer-Verlag, New York}, 1984.

[8] {\sc N. N. Vakhania}, {\sc V. I. Tarieladze} and {\sc S. A. Chobanian},
Probability distributions in Banach spaces,
{\it Nauka, Moscow}, 1985 (in Russian).

 [9] {\sc V. M. Kruglov} and {\sc V. Yu. Korolev},
Limit theorems for random sums,
{\it Moscow University Press, Moscow},
1990 (in Russian).

\end{document}